\documentclass[12pt]{amsart}
\usepackage{amsthm}
\usepackage{amsmath}
\usepackage{amsfonts}
\newtheorem{theorem}{Theorem}[section]
\newtheorem{lemma}[theorem]{Lemma}
\newtheorem{proposition}[theorem]{Proposition}
\newtheorem{corollary}[theorem]{Corollary}
\newtheorem{conjecture}[theorem]{Conjecture}
\newtheorem{definition}[theorem]{Definition}
\newtheorem{example}[theorem]{Example}

\newenvironment{pfmt}{\noindent{\bf Proof of the Main Theorem.}}{}
\newenvironment{ac}{\noindent{\bf Acknowledgements.}}{}

\def\KK{{\mathbb K}}

\def\ZZ{{\mathbb Z}}

\begin{document}
\title{On the structure of $\mu$-classes}
\author{Carlos D' Andrea}
\address{Departamento de Matem\'atica, FCEyN, Universidad de Buenos Aires.
Pabell\'on I, Ciudad Universitaria, 1428 Buenos Aires. Argentina.}
\email{cdandrea@dm.uba.ar}

\begin{abstract}
We prove that, if $\mu<\lfloor\frac n2\rfloor,$ then every rational parametrization of degree $n$ and class $\mu$
is a limit of parametrizations of the same degree and class $\mu+1.$
This property was conjectured in
\cite{CSC}, and its validity allows an explicit description of the variety of parametrizations of degree $n$ and class $\mu,$ for all $(n,\mu).$
\end{abstract}

\maketitle

\section{introduction}
Let $\KK$ be an algebraically closed field. Consider three polynomials
$a=a(t),b=b(t),c=c(t)\in\KK[t],$
of respective degrees $n_a,n_b,n_c,$ such that
$\gcd(a,b,c)=1.$
If $c\neq0,$ then
\begin{equation}\label{pc}
x=\frac{a(t)}{c(t)},\ y=\frac{b(t)}{c(t)}
\end{equation}
is the parametrization of a rational function. Parametrizations of
this form play an important role in computer-aided geometric design
(see for instance \cite[Chapter 6]{CLO}).
\par
A basic object in this area is the syzygy module of the triple $(a,b,c):$
$$Syz(a,b,c):=\{(A,B,C)\in\KK[t]^3:\, Aa+Bb+Cc=0\}.$$
The degree of a vector of polynomials of the form
$(A,B,C)$ is defined as $\max\{\deg(A),
\deg(B),\deg(C)\}.$
Let $n=\max\{n_1,n_b,n_c\},$ and set
$$\mu(a,b,c):=\min\{\deg(A,B,C),\,
(A,B,C)\in Syz(a,b,c)\setminus\{0\}\}.$$ This number is
called the \textit{class} of the triple $(a,b,c).$
It may be regarded as an invariant of the parametric curve
(\ref{pc}), which
provides information about how to implicitize it (see \cite{CS,CZS,CW,CLO,CSC,ZS}).
It is well known (see \cite{CLO,CSC}) that
$\mu\leq[\frac{n}{2}].$ Moreover, the upper bound holds for generic curves of
degree $n.$
\par
In a more general setting, for a given positive integer $n,$
let $\KK[t]_n$ denote the vector space of all
polynomials in $t$ of degree $\leq n,$ and let
${\mathcal P}_n\subset\KK[t]_n^3$ be the subset of $(a,b,c)\in\KK[t]_n^3$
such that $g\neq0,\, \gcd(a,b,c)=1$ and $n=\deg(a,b,c).$
${\mathcal P}_n$ may be regarded as the set of all parametric equations of
planar rational curves or degree $n.$
\par
For every $\mu\in\ZZ_{\geq0},$ set
$${\mathcal P}^\mu_n:=\{(a,b,c)\in{\mathcal P}_n:\,\mu(a,b,c)=\mu\}.$$
It turns out that ${\mathcal P}_n={\mathcal P}^0_n\cup\dots\cup{\mathcal P}^{
\lfloor\frac n2\rfloor}_n.$ In \cite[Theorem 1]{CSC}, it is shown that, for
every $0\leq\mu\leq\lfloor\frac n2\rfloor,$
 the Zariski closure of ${\mathcal P}^\mu_n,$ denoted
$\overline{{\mathcal P}^\mu_n},$
is an irreducible variety of dimension $2n+2\mu+4$ if $\mu<\lfloor\frac n2\rfloor,$ or $3n+3$ if $\mu=\lfloor\frac n2\rfloor.$
\par
In that paper, it is shown that $\overline{{\mathcal P}^{\lfloor\frac n2\rfloor}_n}={\mathcal P}^0_n\cup\dots\cup{\mathcal P}^{
\lfloor\frac n2\rfloor},$ and it was conjectured that, for every $\mu\leq
\lfloor\frac n2\rfloor,$
\begin{equation}\label{fib}
\overline{{\mathcal P}^{\mu}_n}={\mathcal P}^0_n\cup\dots\cup{\mathcal P}^{
\mu}.
\end{equation}
As stated in \cite{CSC}, due to the irreducibility of ${\mathcal P}^\mu_n,$
this is equivalent to the following:
\begin{conjecture}\label{conj}
If $\mu<\lfloor \frac n2\rfloor,$ then every parametrization of class $\mu$ is a limit of parametrizations of class $\mu+1.$
\end{conjecture}
The main result of this paper is a positive answer to this question. This is our main result:
\begin{theorem}\label{mt}
Let $\epsilon$ be a new variable. If $\mu<\lfloor \frac n2\rfloor$ then, for any $(a,b,c)\in{\mathcal P}^\mu_n,$
there exists
$(a_\epsilon,b_\epsilon,c_\epsilon)\in\KK[\epsilon,t]^3$
such that, as univariate polynomials in
$\KK(\epsilon)[t],$ they satisfy the following conditions:
\begin{itemize}
\item  $\gcd(a_\epsilon,b_\epsilon,c_\epsilon)=1;$
\item $n=\max\{\deg(a_\epsilon),\deg(b_\epsilon),\deg(c_\epsilon)\},$
\item $\mu(a_\epsilon,b_\epsilon,c_\epsilon)=\mu(a,b,c)+1,$
\item $(a_\epsilon,b_\epsilon,c_\epsilon)|_{\epsilon=0}= (a,b,c).$
\end{itemize}
\end{theorem}
\begin{definition}
A triple $(a_\epsilon,b_\epsilon,c_\epsilon)\in\KK[\epsilon,t]^3$ satisfying the
hypothesis of Theorem \ref{mt} will be called  an approximating
sequence for $(a,b,c).$
\end{definition}
We may regard approximating sequences as  deformations of $(a,b,c)$ which
``converge'' to the latter when $\epsilon\to0.$
As an immediate result of our construction,
we get that Conjecture \ref{conj} is true, so the equality given in
(\ref{fib}) holds for every $\mu=0,1,
\dots,\lfloor\frac n2\rfloor.$
\par
In the following section we present some auxiliary results, which will
allow us to put a given triple $(a,b,c)\in{\mathcal P}^\mu_n$ in a ``generic
situation''.
In Section \ref{tres}, we will construct the approximating sequence for this generic situation and prove Theorem \ref{mt}.

\section{Auxiliary Results}
In this section we will show that, in order to have some genericity,
we may replace $(a,b,c)$ with $(a+\lambda b,b,c),\,\lambda\in\KK.$ This change ``preserves'' approximating sequences.

The following properties will be useful in the sequel:
\begin{theorem}[\cite{CSC}]\label{cox}
There exists unique polynomials $p_x,p_y,p_w,q_x,q_y,q_w\in\KK[t]$ such that
\begin{enumerate}
\item $\deg(p_x,p_y,p_w)=\mu,\,\deg(q_x,q_y,q_w)=n-\mu;$
\item $a=p_yq_w-p_wq_y,\,b=p_wq_x-p_xq_w,\,c=p_xq_y-p_yq_x;$
\item Every syzygy $(A,B,C)\in Syz(a,b,c)$ can be written uniquely in the
form
$$(A,B,C)=h_1(p_x,p_y,p_w)+h_2(q_x,q_y,q_w),$$
with $h_1,h_2\in\KK[t]$ with $\deg(h_1)\leq\deg(A,B,C)-\mu$ and
$\deg(h_2)\leq\deg(A,B,C)+\mu-n.$
\end{enumerate}
\end{theorem}

\begin{corollary}
$\gcd(p_x,p_y,p_w)=1.$
\end{corollary}
\begin{proof}
As $\frac{(p_x,p_y,p_w)}{\gcd(p_x,p_y,p_w)}\in Syz(a,b,c),$
it turns out that, if the degree of the $\gcd$ is positive, we would have
a syzygy of degree strictly less than $\mu(a,b,c)$ which is a contradiction.
\end{proof}
\begin{lemma}\label{aux3}
For all but finitely many $\lambda\in\KK,\,\gcd(p_y-\lambda p_x,p_w)=1.$
\end{lemma}
\begin{proof}
As the number of roots of $p_w$ is finite, if the claim does
not hold, then there exists $t_0\in\KK$ such that
$p_w(t_0)=0$ and $p_y(t_0)-\lambda p_x(t_0)=0$ for infinitely many
values of $\lambda.$ This implies that $p_y(t_0)=p_x(t_0)=0,$ which is a contradiction with the
corollary.
\end{proof}
\begin{lemma}\label{aux1}
For every $\lambda\in\KK,\,\mu(a,b,c)=\mu(a+\lambda b,b,c).$
\end{lemma}
\begin{proof}
It is straightforward to verify that the morphism $\KK[t]^3\to\KK[t]^3$ given
by $(A,B,C,)\mapsto(A+\lambda B, B, C)$ gives an isomorphism of syzygy modules
which preserves the degree filtration.
To be more precise,
for a fixed $\lambda,$ if we denote with $p^\lambda_x,
p^\lambda_y,\dots,q^\lambda_w$ the polynomials generating the syzygies of
$(a+\lambda b, b, c),$ then
it is straightforward to check that
$$p^\lambda_x=p_x,\, p^\lambda_w=p_w,\,q^\lambda_x=q_x,\,q^\lambda_w=q_w,\,
p^\lambda_y=p_y-\lambda p_x,\, q^\lambda_y=a_y-\lambda q_x.$$
So, we have that $\deg(p^\lambda_x,p^\lambda_y,p^\lambda_w)=\mu(a,b,c).$
Due to Theorem \ref{cox}, this value characterizes $\mu(a+\lambda b,
b,c).$
\end{proof}

\medskip
\begin{proposition}\label{aux2}
Let $\lambda\in\KK.$ If there exists an approximating sequence for
$(a+\lambda b, b, c),$ then there exists an approximating sequence for
$(a,b,c).$
\end{proposition}
\begin{proof}
Let $(\tilde{a_\epsilon},\tilde{b_\epsilon},\tilde{c_\epsilon})\in
\KK[\epsilon,t]^3,$ be an
approximating sequence for $(a+\lambda b, b,c)$
We will show that, if we define
$$a_\epsilon:=\tilde{a_\epsilon}-\lambda\tilde{b_\epsilon}, \ b_\epsilon:=\tilde{b_\epsilon}, \
c_\epsilon:=\tilde{c_\epsilon},$$
we get an approximating sequence for $(a,b,c).$
From the properties of $\tilde{a_\epsilon},\tilde{b_\epsilon},\tilde{c_\epsilon},$
it is easy to check that
$$\gcd(a_\epsilon,b_\epsilon,c_\epsilon)=1,\ n=\max\{\deg(a_\epsilon),\deg(b_\epsilon),\deg(c_\epsilon)\}.$$
Due to Lemma \ref{aux1}, we have that
$$\mu(a_\epsilon,b_\epsilon,c_\epsilon)=\mu(\tilde{a_\epsilon},\tilde{b_\epsilon},\tilde{c_\epsilon})=
\mu(a+\lambda b,b,c)=\mu(a,b,c)+1,$$
and it is straightforward to see that
$(a_\epsilon,b_\epsilon,c_\epsilon)|_{\epsilon=0}=(a,b,c).$
\end{proof}

\section{Construction of the Approximating Sequence}
\label{tres}
In this section we will fix $(a,b,c)\in{\mathcal P}^\mu_n$ such that
$\mu+1\leq\lfloor\frac{n}{2}\rfloor,$ and will construct an approximating
sequence for $(a+\lambda b, b,c)$ for an appropiate $\lambda\in\KK.$ Proposition
\ref{aux2} implies that there exists an approximating sequence for $(a,b,c).$
\par
Suppose w.l.o.g. that $\deg(a)$ is positive. We will consider a family
$(a_\epsilon,b_\epsilon,c_\epsilon)$ of the form
\begin{equation}\label{aps}
a_\epsilon:=a(1+\frac{\epsilon}{t-\alpha}),\,b_\epsilon:=b,\,
c_\epsilon:=c,
\end{equation}
where $\alpha\in\KK$ is a root of $a$ having properties to be
described shortly.
It is clear that, as $\epsilon$ goes to zero, the deformed family converges to
$(a,b,c).$
Moreover, it is straightforward to check that
$n(a_\epsilon,b_\epsilon,c_\epsilon)=n(a,b,c)$ and that
$\gcd(a_\epsilon,b_\epsilon,c_\epsilon)=\gcd(a,b,c)=1.$
In order to have $\mu(a_\epsilon,b_\epsilon,c_\epsilon)=\mu(a,b,c)+1,$
we choose $\alpha$ such that $a(\alpha)=0$ and either $p_y(\alpha)\neq0$
or $p_w(\alpha)\neq0.$
\par It is not true that the latter condition can always be acomplished as
the following cautionary example shows.
\begin{example}\label{exastuto}{\rm
Set $p_x=q_w=t,\,p_y=p_w=t-1,\,q_x=t^4,q_y=1.$ Then,
$a=(t-1)^2,\,b=t^5-t^4-t^2,\,c=-t^5+t^4+t,$ and we have that there is a unique
root of $a$ which is also the same root of $p_y$ and $p_w.$}
\end{example}
If we are in the situation that every root of $a$ is also a root of $p_y$ and
$p_w,$ then we will construct an approximating sequence for
$(a+\lambda b,b,c)$ for $\lambda\in\KK^*$ such that
$\deg(a+\lambda b)>0.$
\begin{lemma}\label{aux4}
There is $\lambda\in\KK$ such that $\deg(a+\lambda b)>0,$ and $a+\lambda b$ has a root $\alpha$
such that either $p_y^\lambda(\alpha)\neq0$
or $p_w^\lambda(\alpha)\neq0.$
\end{lemma}
\begin{proof}
If every root of $a$ is a root of $p_y$ and $p_w,$ then  pick one of those $\lambda$ satisfying Lemma \ref{aux3},
 and denote, as in the proof of Proposition
\ref{aux3}, with $p^\lambda_x,
p^\lambda_y,\dots,q^\lambda_w$ the polynomials generating the syzygies of
$(a+\lambda b, b, c),$ we get that $p^\lambda_y=p_y-\lambda p_x,$ and
$p^\lambda_w=p_w.$
\par
As $\gcd(p^\lambda_y,p^\lambda_w)=\gcd(p_y-\lambda p_x,p_w)=1$ by Lemma
\ref{aux3}, then any root
of $a+\lambda b$ may be chosen provided that $\deg(a+\lambda b)$ is positive.
\par
In order to see the latter,
let $\lambda_1\ne\lambda_2$ both satisfying Lemma \ref{aux3}. If
$a+\lambda_1 b$ and $a+\lambda_2 b$ were both constant, then so would $a$ and
$b.$
\end{proof}

\medskip
\begin{pfmt}
We may suppose then, w.l.o.g. that there exists $\alpha\in\KK$
such that is a root of $a$ but is not a common root
of $p_y$ and $p_w.$
\par
Due to the remarks made at the beginning of this section, it only remains to
prove that $\mu(a_\epsilon,b_\epsilon,c_\epsilon)=\mu+1.$
First, observe that
$$p_x(t-\alpha)a_\epsilon+p_y(t-\alpha+\epsilon)b_\epsilon+p_w(t-\alpha+
\epsilon)c_\epsilon=(t-\alpha+\epsilon)(p_xa+p_yb+p_wc)=0,
$$
so $((t-\alpha)p_x,(t-\alpha+\epsilon)p_y,(t-\alpha+\epsilon)p_w)$ is a
nontrivial syzygy on $(a_\epsilon,b_\epsilon,c_\epsilon)$ and we deduce
that $\mu(a_\epsilon,b_\epsilon,c_\epsilon)\leq \mu(a,b,c)+1.$
\par
Suppose that the inequality holds strictly. Then, there must be $(A,B,C)$
of degree bounded by $\mu$ such that $Aa_\epsilon+Bb_\epsilon+Cc_\epsilon=0.$
As
$$\begin{array}{l}
A(t-\alpha+\epsilon)a+B(t-\alpha)b+C(t-\alpha)c=\\
A(t-\alpha)\frac{t-\alpha+\epsilon}{t-\alpha}a+B(t-\alpha)b+C(t-\alpha)c=\\
A(t-\alpha)a_\epsilon+B(t-\alpha)b_\epsilon+C(t-\alpha)c_\epsilon=
(t-\alpha)(Aa_\epsilon+Bb_\epsilon+Cc_\epsilon)=0,
\end{array}$$
we have that $(A(t-\alpha+\epsilon),B(t-\alpha),C(t-\alpha))\in Syz(a,b,c),$
and it is a nontrivial syzygy of degree bounded by $\mu+1.$
So, due to the last item of Theorem \ref{cox}, there exists $h_1,h_2
\in\KK[t]$ such that
$$(A(t-\alpha+\epsilon),B(t-\alpha),C(t-\alpha))=
h_1(p_x,p_y,p_w)+h_2(q_x,q_y,q_w),$$
with $\deg(h_2)\leq (\mu+1)+\mu-n.$
As we have supposed $\mu+1\leq\frac n2,$ then 
$2\mu+1-n\leq-1,$ so $h_2$ must be identically zero. Then, we have that
$(A(t-\alpha+\epsilon),B(t-\alpha),C(t-\alpha))=h_1(t)(p_x,p_y,p_w),$
with $\deg(h_1)=1.$
\par Comparing the first coordinate of the last equality, we have that
$h_1(t)=\lambda'(t-\alpha+\epsilon),\,\lambda'\in\KK\setminus\{0\}.$
In addition, if we compare the two last coordinates, we will have that
$(t-\alpha)$ divides both $p_y$ and $p_w,$ which contradicts our hypothesis.
So, we have that $\mu(a_\epsilon,b_\epsilon,c_\epsilon)=\mu(a,b,c)+1$ and
the claim holds.
\end{pfmt}

\begin{example}{\rm
We will show here that, in order to construct to approximation family, if we take any root $\alpha$ of $a$ without
imposing the condition that it should not be simultaneously
root of $p_y$ and $p_w,$ the sequence (\ref{aps}) may fail to have degree
$\mu+1.$ Consider again, as in Example \ref{exastuto},
$$
a:=(t-1)^2,\,b=t^5-t^4-t^2,\,c=-t^5+t^4+t.$$
If we take as $\alpha$ the unique root of $a,$ then we will have that
$$a_\epsilon=(t-1)^2+\epsilon(t-1),\,b_\epsilon=b,\,c_\epsilon=c.$$
This is not an approximating sequence, because $\mu(a,b,c)=1$ and,
as
$$ta_\epsilon+(t-1+\epsilon)b_\epsilon+(t-1+\epsilon)c_\epsilon=0,$$
then we have that $\mu(a_\epsilon,b_\epsilon,c_\epsilon)=1.$
\par
In this case, we may choose $\lambda=-\frac{1}{12}.$ One can check that
$\alpha=2$ is a root of $a+\lambda b,$ and with this choice of $\alpha,$
we get the following family:
$$\begin{array}{l}
a_\epsilon=-\frac{t^5}{12}+\frac{t^4}{12}+\frac{13}{12}t^2-2t+1+\epsilon(
-\frac{1}{12}t^4-\frac{1}{12}t^3-\frac16t^2+\frac34t-\frac12),\\
b_\epsilon=b,\,c_\epsilon=c.
\end{array}$$
An easy calculation shows that this deformed family has class
equal to two, so we get an approximating family of $(a+\frac12b,b,c),$ and
hence
$(a+\epsilon(-\frac{1}{12}t^4-\frac{1}{12}t^3-\frac16t^2+\frac34t-\frac12),b,c)$ is
an approximating family of $(a,b,c).$}
\end{example}
\begin{ac}
I would like to thank David Cox for his careful reading of the manuscript and many helpful suggestions.
This research was conducted while I was
a postdoctoral fellow at the
Institut National de Recherche en
Informatique et en Automatique (INRIA) in
Sophia-Antipolis, France,  partially supported by
Action A00E02 of the ECOS-SeTCIP French-Argentina bilateral collaboration.
\end{ac}

\end{document}